\documentclass[]{IEEEtran}

\usepackage{graphicx}
\usepackage{amsmath}
\usepackage{amssymb}

\usepackage[usenames,dvipsnames]{xcolor}
\usepackage{dsfont}
\usepackage{balance}
\usepackage[usenames,dvipsnames]{xcolor}
\makeatletter
\def\amsbb{\use@mathgroup \M@U \symAMSb}
\makeatother

\usepackage{bbold}

\newtheorem{thm}{Theorem}
\newtheorem{cor}{Corollary}
\newtheorem{lem}{Lemma}
\newtheorem{rem}{Remark}
\newtheorem{defn}{Definition}

\title{Delay-independent stability in monotone systems}

\author{Eoin~Devane and Ioannis~Lestas
\thanks{This work was supported in part by the Engineering and Physical Sciences Research Council (EPSRC). Paper~\cite{acc2015} is a preliminary conference version of this manuscript. This manuscript includes detailed proofs of Theorem~\ref{theo2}, Lemma~\ref{theo1}, {\color{black}and Corollaries~\ref{cor1} and~\ref{cor2}} which are omitted in {\color{black}the final version of}~\cite{acc2015}, and also additional examples and discussion.}
\thanks{Eoin Devane is with the {\color{black}Department of Pure Mathematics and Mathematical Statistics}, University of Cambridge, Wilberforce Road, Cambridge, CB3 0WA, United Kingdom; \tt\small esmd2@cam.ac.uk}
\thanks{Ioannis Lestas is with the Department of Engineering, University of Cambridge, Trumpington Street, Cambridge, CB2 1PZ, United Kingdom; \tt\small icl20@cam.ac.uk}}

\begin{document}

\maketitle

\begin{abstract}
Monotone systems comprise an important class of dynamical systems that are of interest both for their wide applicability and because of their interesting mathematical properties. It is known that under the property of quasimonotonicity time-delayed systems become monotone, and some remarkable properties have been reported for such systems. These include, for example, the fact that for linear systems global asymptotic stability of the undelayed system implies global asymptotic stability for the delayed system under arbitrary bounded delays. Nevertheless, extensions to nonlinear systems have thus far relied {\color{black}primarily on the conditions of homogeneity and subhomogeneity}, and it has been conjectured that these can be relaxed. Our aim in this paper is to show that this is feasible for a general class of nonlinear monotone systems by deriving convergence results in which simple properties of the undelayed system lead to delay-independent stability. In particular, one of our results shows that if the undelayed system has a convergent trajectory that is unbounded {\color{black}in all components} as $t\rightarrow-\infty$, then the system is globally asymptotically stable for arbitrary {\color{black}bounded} time-varying delays. This follows from a more general result derived in the paper that allows to quantify delay-independent regions of attraction, which can be used to  prove global asymptotic stability for various classes of systems. These also recover various known delay-independent stability results that are discussed within the paper.
\end{abstract}


\section{Introduction} \label{sec:intro}
The class of monotone systems, which comprises dynamical systems that preserve an order on the state space, is of significant interest both due to the strong theoretical properties of such systems and for their suitability for modeling numerous physical systems. Areas in which monotonicity properties have frequently been exploited include population dynamics~\cite{kolmogorov}, consensus protocols~\cite{altafini}, and communication systems~\cite{yates}, emphasizing the breadth of applicability of {\color{black}monotone} systems. The seminal papers of Hirsch, beginning with~\cite{hirsch}, established a powerful general theory for monotone systems, demonstrating that the order-preservation property endows such systems with a rich mathematical structure. A thorough review of this theory can be found in~\cite{hirschsmith}, which details in particular the powerful results of generic convergence that can be obtained.

It is well-known that delay differential equations become monotone under the property of quasimonotonicity~\cite{smith}. In the modeling of many physical problems, time-delays play an important role and can have a significant effect on the system behavior. In particular, time-delays can often destabilize a system and prevent convergence from occurring. It is therefore an important problem to ascertain how the presence of time-delays might affect the behavior of systems that are asymptotically stable in the undelayed setting. Quasimonotonicity alone is not in general sufficient to be able to guarantee global asymptotic stability. Furthermore, despite the monotonicity property, direct comparisons between the trajectories of the undelayed and the delayed system are in general nontrivial to establish, as the presence of delays can significantly change the behavior of trajectories. Therefore, there have been a number of attempts to resolve this problem for particular types of systems, through use of monotonicity combined with other system properties. For linear systems there exist strong results, many of which are based upon versions of the Perron--Frobenius Theorem, guaranteeing delay-independent global asymptotic stability and even exponential convergence~rates~\cite{farina,liu}. In recent interesting works, these results have been extended to nonlinear systems satisfying properties of homogeneity~\cite{mason,bokharaie} and subhomogeneity~\cite{bokharaie2}. Furthermore, it was shown in~\cite{bokharaie2} that the subhomogeneity condition can be sufficiently strong to imply delay-independent stability even without monotonicity. However, since homogeneity is closely linked to linearity, it is desirable to provide extensions, as was done in the planar case in~\cite{bokharaiecorrec}, to a more general nonlinear monotone setting.

Within this paper we investigate what conclusions can be drawn about the behavior of the delayed system using the quasimonotonicity property alone, together with simple properties of the undelayed system. In particular, one of our results is to show that whenever the system without any delay admits a solution that is convergent and is unbounded in all components in negative time, then global asymptotic stability is guaranteed for arbitrary bounded delays, which can be heterogeneous and time-varying. This follows from a more general result derived in the paper, which explicitly constructs delay-independent regions of attraction in terms of {\color{black}convergent} trajectories of the corresponding undelayed system. The significance of these regions is that they allow us to generate proofs of global asymptotic stability for various classes of systems, {\color{black}recovering also known delay-independent stability results for linear/homogeneous/subhomogeneous systems and convergence for systems with vector fields that are negative at a prescribed point.}

The paper is structured as follows. The notation that will be used within the paper is described in Section~\ref{sec:notation}. The problem formulation is given in Section~\ref{sec:probform}. The main results are then stated in Section~\ref{sec:convres} and their proofs are given in Section~\ref{sec:proofs}. {\color{black}Section~\ref{sec:examples} includes a discussion on the significance of the results in the paper, and conclusions are drawn in Section~\ref{sec:conclusions}.} {\color{black}In Appendix~A we review some classical results on monotone systems that are used within the paper, and in Appendix~B we derive some extensions of our results to nonpositive monotone systems.}

\section{Notation} \label{sec:notation}

The set of $n$-dimensional vectors with real entries is denoted by~$\mathds{R}^n$, and
inequalities in $\mathds{R}^n$ are defined as follows: $x \ge y$ means $x_i \ge y_i$ {\color{black}for all $i$}; $x > y$ means $x_i \ge y_i$ {\color{black}for all $i$} and $x \neq y$;  $x \gg y$ means $x_i > y_i$ {\color{black}for all $i$}. {\color{black}The nonnegative orthant $\{x \in \mathds{R}^n \colon x \ge 0\}$ is denoted by $\mathds{R}_+^n$.}

We use $\mathcal{C}([-r,0],\mathds{R}^n)$ to denote the Banach space of continuous functions mapping the interval $[-r,0]$ into $\mathds{R}^n$. Inequalities in $\mathcal{C}([-r,0],\mathds{R}^n)$ are treated pointwise, for example $\phi \ge \psi$ means $\phi(\theta) \ge \psi(\theta)$ for all $\theta \in [-r,0]$. The function $\mathds{1}$ represents the constant mapping $\mathds{1}: [-r,0] \mapsto 1_n$, where $1_n \in \mathds{R}^n$ has all entries equal to 1. {\color{black}The product $v \mathds{1}$\color{black}, where $v\in\mathds{R}^n$, denotes the constant mapping $v \mathds{1}: [-r,0] \mapsto v$.} We let $x_t$ denote the function in $\mathcal{C}([-r,0],\mathds{R}^n)$ given by $x_t(\theta) = x(t+\theta)$ for all $\theta \in [-r,0]$ for a solution $x(t)$ of the delayed system that will be introduced in~\eqref{delsys}. The subscript in $x_k$ will sometimes also be used to denote the $k^\text{th}$ component of the variable $x(t)$, but this notational ambiguity will be clear from the context.

Finally, the expression $\frac{dx}{dt}(\tau)$ will be used to denote the evaluation of the time-derivative of $x$ at~$t = \tau$,~$\frac{d}{dt}(x(t))\big|_{t = \tau}$.

\section{Problem formulation} \label{sec:probform}
Throughout the paper we investigate the stability properties of a general class of monotone systems with time-varying {\color{black}delays,} which we define in this section.

We consider a function $f: \mathds{R}^n \to \mathds{R}^n$ and a functional $g^t: \mathcal{C}([-r,0],\mathds{R}^n) \to \mathds{R}^n$ that satisfy the {\color{black}conditions}
\begin{align}
&x \le y \text{ and } x_i = y_i \Rightarrow f_i(x) \le f_i(y), \label{fmon} \\
&\phi \le \psi \Rightarrow g^t(\phi) \le g^t(\psi), \label{gmon}
\end{align}
for all $x, y \in \mathds{R}^n$, all $\phi, \psi \in \mathcal{C}([-r,0],\mathds{R}^n)$, and all $t~\ge~0$. We then formulate the {\color{black}delayed system
\begin{equation}
\frac{dx}{dt}(t) = f(x(t)) + g^t(x_t), \label{delsys}
\end{equation}
which will be analyzed throughout the paper. We also {\color{black}consider the} {\color{black}corresponding undelayed system}
\begin{equation}
\frac{dx}{dt}(t) = f(x(t)) + g^t(x(t) \mathds{1}). \label{undelsys}
\end{equation}}

{\color{black}We} assume that the time-dependence $t \mapsto g^t(\phi)$ is continuous for all fixed $\phi \in \mathcal{C}([-r,0],\mathds{R}^n)$ and that the mappings $x \mapsto f(x)$, $\phi \mapsto g^t(\phi)$ are both locally Lipschitz for all fixed $t \ge 0$. The systems~\eqref{delsys} and~\eqref{undelsys} are then guaranteed to have unique solutions throughout $\mathds{R}_+ \times \mathds{R}^n$ which exhibit continuous dependence on the initial data~\cite[Theorems 2.2.1, 2.2.2, and 2.2.3]{hale}.

{\color{black}It is well known that conditions~\eqref{fmon} and~\eqref{gmon} render system\footnote{\color{black}Note that system~\eqref{undelsys} is also monotone as this is a special case of \eqref{delsys}.}~\eqref{delsys} monotone, i.e. if two initial conditions {\color{black}$x_{t_0}, \tilde{x}_{t_0} \in \mathcal{C}([-r,0],\mathds{R}^n)$ satisfy $x_{t_0} \le \tilde{x}_{t_0}$, then the corresponding solutions of~\eqref{delsys} satisfy $x(t) \le \tilde{x}(t)$} for all $t \ge t_0$~\cite{smith}. This follows from the fact that the vector field in~\eqref{delsys} satisfies the property of quasimonotonicity (see more detailed discussion of this in {\color{black}Appendix A)}.}

{\color{black}Furthermore}, it is assumed that given any $v \in \mathds{R}^n$, the map $t \mapsto g^t(v\mathds{1})$ is constant, independent of time $t$. {\color{black}This is} always satisfied if the time-variation in $g^t$ is due to time-varying delays.~Therefore, system~\eqref{delsys}, whose stability properties will be investigated, is a monotone system that can have {\color{black}\emph{arbitrary}} bounded time-varying~delays. 

The equilibria of the systems~\eqref{delsys} and~\eqref{undelsys} are known to be equivalent~\cite{smith}. The problem we are interested in is to relate the stability properties of an equilibrium of the delayed system~\eqref{delsys}, which is infinite-dimensional, to those of this equilibrium for the finite-dimensional undelayed system~\eqref{undelsys}. For simplicity, we will suppose that the equilibrium under consideration is located at the origin, whence~\eqref{delsys} and~\eqref{undelsys} become positive systems (i.e. the nonnegative orthant $\mathds{R}_+^n$ is positively invariant\footnote{By positive invariance of a subset $X \subseteq \mathds{R}^n$ with respect to~\eqref{delsys}, we mean that $x_{t_0} \in \mathcal{C}([-r, 0], X)$ implies $x(t) \in X$ for all $t \ge t_0$. Likewise, when we say that $Y \subseteq \mathds{R}^n$ is a region of attraction of an equilibrium $x^*$ for~\eqref{delsys}, we mean that $x_{t_0} \in \mathcal{C}([-r, 0], Y)$ implies $\lim_{t\to\infty} x(t) = x^*$.} with respect to both systems) due to the monotonicity properties~\eqref{fmon} and~\eqref{gmon}. As is common in the analysis of positive systems, we use standard definitions for stability (e.g. Definition 5.1.1 in~\cite{hale}) but with any region for $x_t$ used in the definition restricted to the cone $\mathcal{C}([-r, 0], \mathds{R}_+^n)$. {\color{black}The origin} is defined to be a {\em globally asymptotically stable} equilibrium of~\eqref{delsys} if it is stable and $\lim_{t\to\infty} x(t) = 0$ for any solution $x(t)$ of~\eqref{delsys} with initial condition $x_{t_0}$ in $\mathcal{C}([-r, 0], \mathds{R}_+^n)$. {\color{black}In Appendix B we show how the stability results in Section~\ref{sec:convres} can then be extended to establish convergence results also {\color{black}for monotone systems} that are not positive.}

\section{Convergence results} \label{sec:convres}

In this section we will state the main results of the paper. The proofs of these results are given in Section~\ref{sec:proofs} and further discussion on their significance is included in Section~\ref{sec:examples}.

One of our main results is given in the following theorem, which states that the origin in~\eqref{delsys} is globally asymptotically stable if there exists any single convergent trajectory of the undelayed system~\eqref{undelsys} that is unbounded in all components in negative time.

\begin{thm} \label{theo2}
Suppose that the undelayed system~\eqref{undelsys} admits a solution $y(t)$ satisfying $\lim_{t \to -\infty} y_i(t) = \infty$ and $\lim_{t \to \infty} y_i(t) = 0$ for all components $i$. Then the origin is globally asymptotically stable for the delayed system~\eqref{delsys}.
\end{thm}
\begin{rem} \normalfont
It should be noted that system~\eqref{delsys} is allowed to have arbitrary bounded delays, which can be heterogeneous and time-varying.

\end{rem}
\begin{rem} \normalfont
Theorem~\ref{theo2} can be seen as addressing the problem of when global asymptotic stability of an undelayed monotone system also implies delay-independent global asymptotic stability. As discussed in the introduction, existing results in the literature have relied primarily on the properties of homogeneity and subhomogeneity, and it has been conjectured in~\cite{bokharaie} that these assumptions can potentially be relaxed. Theorem~\ref{theo2} shows that this is possible for systems having a backward-time divergent trajectory.
\end{rem}

\begin{rem} \normalfont
The conditions on $y(t)$ given in Theorem~\ref{theo2} trivially imply global asymptotic stability for an undelayed monotone system. Nevertheless, proving that stability is not compromised under arbitrary bounded delays is a more involved problem. This is due to the fact that direct comparisons between the trajectories of the undelayed and the corresponding delayed system are in general nontrivial to establish. A main part of the proof of the theorem is therefore associated with the derivation of results that allow us to quantify such connections.
\end{rem}

Theorem~\ref{theo2} follows from Lemma~\ref{theo1} stated below, which is a more general technical result that allows the explicit deduction of delay-independent regions of attraction for~\eqref{delsys} from convergent trajectories of the undelayed system~\eqref{undelsys}.

\begin{lem} \label{theo1}
Suppose that the undelayed system~\eqref{undelsys} admits a solution $y(t)$ satisfying $y_i(0) > 0$ and $\lim_{t \to \infty} y_i(t) = 0$ for all components~$i$. Let $\zeta \ge 0$ be a point defined in terms of $y(t)$ by means of equations~\eqref{Tdef}--\eqref{zdef}. Then any solution of the delayed system~\eqref{delsys} with initial condition $x_{t_0} \in \mathcal{C}([-r,0],[0,\zeta])$ satisfies $\lim_{t \to \infty} x(t)~=~0$.
\end{lem}

\begin{rem} \normalfont
The significance of the point $\zeta$ will become clearer in the proof of Lemma~\ref{theo1} and in the discussion in Section~\ref{sec:examples}. In particular, a main significance {\color{black}of Lemma~\ref{theo1} is that for many classes of systems the point $\zeta$ can be chosen to be arbitrarily large, thereby guaranteeing delay-independent global asymptotic stability. This leads to the proof of Theorem~\ref{theo2} and it also recovers known delay-independent stability results for linear/homogeneous/subhomogeneous systems.}
\end{rem}

Lemma~\ref{theo1} has the following immediate corollary for systems with continuously differentiable vector fields. In particular, it is guaranteed in this case that $\zeta\gg0$, which ensures delay-independent asymptotic stability of the system~\eqref{delsys}. Note that this condition on the vector field is not necessary to deduce asymptotic stability and is also not an assumption in Theorem~\ref{theo2}.

\begin{cor} \label{cor1}
Suppose that $f$ and $g^t$ are both continuously differentiable for all fixed $t$ and that the undelayed system~\eqref{undelsys} admits a solution $y(t)$ satisfying $y_i(0) > 0$ and $\lim_{t \to \infty} y_i(t) = 0$ for all {\color{black}components} $i$. Then the point $\zeta$ in Lemma~\ref{theo1} can be chosen to be strictly positive in all components, and the origin is asymptotically stable for the delayed system~\eqref{delsys}.
\end{cor}

Corollary~\ref{cor1} implies also the following simple corollary.
\begin{cor} \label{cor2}
Suppose that $f$ and $g^t$ are both continuously differentiable for all fixed $t$ and that the origin is asymptotically stable for the undelayed system~\eqref{undelsys}. Then the origin is also asymptotically stable for the delayed system~\eqref{delsys}.
\end{cor}

\section{Proofs of the results} \label{sec:proofs}

We now present the proofs of the main results stated in Section~\ref{sec:convres}. We begin with the derivation of Lemma~\ref{theo1} and then use this to prove Corollaries~\ref{cor1} {\color{black}and \ref{cor2}} and Theorem~\ref{theo2}.

Within the proof of Lemma~\ref{theo1}, we will make frequent use of two important results from~\cite{smith}, namely Theorem 5.1.1 and Corollary 5.2.2. These results are quoted in the appendix as Lemmas~\ref{st511} and~\ref{sc522}.

\begin{IEEEproof}[Proof of Lemma~\ref{theo1}]
We first define, in terms of the known trajectory $y(t)$ of~\eqref{undelsys}, various quantities that will be used throughout the proof of the lemma, and which are also needed in order to define the point $\zeta$.

\begin{itemize}
\item Firstly, for all $x \le y(0)$, we define, as in~\cite{rantzer}, the functions
\begin{equation} \label{Tdef}
\begin{aligned}
&T_i(x_i) := \sup \{\tau \colon x_i \le y_i(s) \; \forall s \in [0, \tau]\},\\
&V_i(x_i) := e^{-T_i(x)}, \\
&T(x) := \min_i T_i(x_i),\ \text{and} \ V(x) := \max_i V_i(x_i)
\end{aligned}
\end{equation}
in terms of the solution $y(t)$ of~\eqref{undelsys}. Then $V(x) = e^{-T(x)}$ and, according to the given asymptotic properties of $y(t)$, $V : [0,y(0)] \to \mathds{R}_+$ is well-defined. Moreover, due to the continuity of trajectories of~\eqref{undelsys},
each $T_i$, and hence $V_i$, is left-continuous.

\item We also introduce the mapping $h_i : \mathds{R}_+ \to \mathds{R}_+$,
\begin{equation} \label{hdef}
\begin{aligned}
h_i(u) := \inf \{ \tau \ge u \colon \exists \, a \in [0,y_i(0)]& \\
&\hspace{-5em}\text{ such that } T_i(a) = \tau \},
\end{aligned}
\end{equation}
which is well-defined by the asymptotic properties~of~$y(t)$.

\item Let $t^p$ denote a time for which
    \begin{equation}\label{eq:tpdef}
    y(t^p) \ll y(0).
    \end{equation}
Using the map $h$ in~\eqref{hdef}, we define the point $\zeta \in \mathds{R}_+^n$ by
\begin{equation}
\zeta_i := y_i(h_i(t^p)) \label{zdef}
\end{equation}
for all $i$. Note\footnote{{\color{black}If the vector field is also continuously differentiable, then it follows from Corollary~\ref{cor1} that $\zeta$ can be chosen such that $\zeta \gg 0$.} It can easily be seen from the definition of $\zeta$ that the weaker condition $y(t) \gg 0 \ \forall t \in [0,t^p]$ is also sufficient for the existence of a $\zeta$ that is strictly positive.}\textsuperscript{,}\footnote{It should be noted that for a given trajectory $y(t)$ there could exist multiple points $\zeta$ (and times $t^p$) that satisfy~\eqref{eq:tpdef},~\eqref{zdef}, with each such point leading to a different estimate of the region of attraction. It will be seen in the proof of Theorem~\ref{theo2}, in Corollary~\ref{homogcompcor}, and in {\color{black}the example in Section~\ref{sec:examplesgas}} that for various classes of systems the point $\zeta$ can be appropriately chosen so as to also deduce global asymptotic stability.} that $\zeta \ge 0$.

\end{itemize}
The significance of these notions will become clear throughout the proof.
We now give a brief sketch of the proof to facilitate its readability. The function $V(x)$ will be used to construct a compact positively invariant set $S$ for the delayed system. We will then show that $\zeta$ is a maximal point in this set. This property will be used to show that the vector field is nonpositive at $\zeta$, for an appropriately constructed bounding monotone system. Hence the desired convergence result can then be deduced.

The proof is now given in several steps:

\textbf{Step 1:} We first show that the level set
\begin{equation}
S := \{x \in [0,y(0)] \colon V(x) \le e^{-t^p}\} \label{Sdef}
\end{equation}
is positively invariant with respect to the delayed system~\eqref{delsys}.

To establish this, we use an argument similar to Lyapunov--Razumikhin analysis to demonstrate that any trajectory initially within $S$ must remain within $S$ for all future time. We begin by proving the following lemma.

\begin{lem}
\label{lem:Razum}
Suppose that $V(x(s+\theta)) \le V(x(s))$ for all $\theta \in [-r,0]$, for a solution $x(t)$ of the delayed system~\eqref{delsys} and a time $s$ for which the segment $x_s \in \mathcal{C}([-r,0],[0,y(0)])$. Let~$k$ be {\color{black}any} component of $x(s)$ such that $V(x(s)) = V_k(x_k(s))$. Then
\begin{equation*}
\frac{dx_k}{dt}(s) \le 0.
\end{equation*}
\end{lem}

\begin{IEEEproof}[Proof of Lemma~\ref{lem:Razum}]
The conditions stated in Lemma~\ref{lem:Razum} imply that $T(x(s+\theta)) \ge T(x(s))$ for all $\theta \in [-r,0]$. This means that whenever $x(s) \le y(\sigma)$ for all $\sigma \in [0,\tilde{s}]$, we must have $x(s+\theta) \le y(\tilde{s})$. In particular, for $\tilde{s} = T(x(s))$ we therefore get $x(s+\theta) \le y(T(x(s)))$ for all $\theta \in [-r,0]$. Now, by definition, $V(x(s)) = V_k(x_k(s))$ implies $y_k(T(x(s))) = x_k(s)$, whence~\eqref{fmon} and~\eqref{gmon} respectively give
\begin{align}
&f_k(x(s)) \le f_k\big(y(T(x(s)))\big), \label{fmon1}\\
&g^s_k(x_s) \le g^s_k\big(y(T(x(s))) \mathds{1}\big). \label{gmon1}
\end{align}
By the particular choice of component $k$ and the definition~\eqref{Tdef}, we have $\tfrac{dy_k}{dt} (T(x(s))) \le 0$, whence {\color{black}the properties}~\eqref{fmon1} and~\eqref{gmon1} and the time-invariance of $g^t$ for constant arguments give
\begin{align*}
\hspace{2.8em}\frac{dx_k}{dt}(s) &\le f_k\big(y(T(x(s)))) + g^s_k\big(y(T(x(s))) \mathds{1}\big) \nonumber \\
&= f_k\big(y(T(x(s)))) + g^{T(x(s))}_k\big(y(T(x(s))) \mathds{1}\big) \nonumber \\
&= \frac{d y_k}{dt} (T(x(s))) \le 0.
\end{align*}\end{IEEEproof}

Let us now consider the set defined in~\eqref{Sdef}. The asymptotic properties of $y(t)$ and the left-continuity of $V$ ensure that $S$ is compact. By this compactness, the only way in which the desired positive invariance could be violated is if there exists a trajectory $x(t)$ of~\eqref{delsys} and a time $\tau$ such that $x(\sigma) \in S$ for all $\sigma \le \tau$ and for which given any $\epsilon > 0$ there exists $\upsilon \in (\tau,\tau+\epsilon)$ with $x(\upsilon) \notin S$. Since $V(x(\tau)) \le e^{-t^p}$ implies $T(x(\tau)) \ge t^p$ and we know that $y(t^p) \ll y(0)$, it follows that $x(\tau) \ll y(0)$. We can thus always choose $\epsilon$ sufficiently small that $x(\sigma) \in [0,y(0)]$, so that $V(x(\sigma))$ is defined, for all $\sigma \in [\tau,\upsilon]$. Therefore it must be the case that $V(x(\sigma)) \le e^{-t^p}$ for all $\sigma \le \tau$ and $V(x(\upsilon)) > e^{-t^p}$. From the continuous differentiability of the trajectory $x(t)$ and the left-continuity and strictly increasing nature of each $V_i$, it follows that there must exist an interval $(\gamma_1,\gamma_2) \subseteq [\tau,\upsilon]$ such that $V(x(\gamma_1)) \ge \max \{e^{-t^p}, \sup_{\sigma\in[\tau,\gamma_1]}V(x(\sigma)) \}$ and on which the mapping $\sigma \mapsto V(x(\sigma))$ is strictly increasing. Moreover, we know from the continuous differentiability of the trajectory $x=y(t)$ of~\eqref{undelsys} that each $V_i$ is continuous everywhere except on a (possibly infinite) discrete set of isolated points. This means that there must exist some $\gamma_3 \in (\gamma_1,\gamma_2)$ such that $V(x(\sigma)) = V_k(x_k(\sigma))$ for all $\sigma \in (\gamma_1,\gamma_3)$ for some fixed choice of $k$. It thus follows that the mapping $\sigma \mapsto x_k(\sigma)$ is strictly increasing on $(\gamma_1,\gamma_3)$, whence the Mean Value Theorem ensures the existence of some $\beta \in (\gamma_1,\gamma_3)$ at which $\dot{x}_k(\beta) > 0$. But, by construction, we know that $V(x(\beta)) \ge V(x(\beta + \theta))$ for all $\theta \in [-r,0]$ and $V(x(\beta)) = V_k(x_k(\beta))$, so this positive time-derivative is in contradiction to Lemma~\ref{lem:Razum}. Therefore, $S$ must be positively invariant with respect to~\eqref{delsys}.

\textbf{Step 2:} We now show that $\zeta$ is a maximal point of the set $S$. In particular, we show that $S \subseteq [0,\zeta]$ and $\zeta \in S$.

To do this, consider any $x \in S$. According to the definition~\eqref{hdef}, it must follow from $T_i(x_i) \ge t^p$ that $T_i(x_i) \ge h_i(t^p)$, whence $x_i \le y_i(h_i(t^p)) = \zeta_i$. Since this holds for all $i$, we immediately conclude that $S \subseteq [0,\zeta]$.

Additionally, since each $T_i$ is a nonincreasing, left-continuous function on a compact set, the image $T_i([0,y_i(0)])$ contains all of its right-limits. Therefore, by~\eqref{hdef}, for any $i$ there must exist some $a \in [0,y_i(0)]$ such that $T_i(a) = h_i(t^p)$. It thus holds by~\eqref{Tdef} that $a = y_i(h_i(t^p)) = \zeta_i$. This shows that $\zeta$ lies within the domain of definition of $T$, $V$ and moreover, because $h_i(t^p) \ge t^p$, that $T(\zeta) \ge t^p$. Therefore, $\zeta \in~S$.

\textbf{Step 3:} We now use the results of Steps 1 and 2 to guarantee that the set $[0,\zeta]$ lies within the region of attraction of the origin of~\eqref{delsys}.

From Steps 1 and 2, we conclude that the sublevel set $S$ satisfies positive invariance with respect to~\eqref{delsys}, is contained within the interval~$[0,\zeta]$, and contains the point $\zeta$. Applying positive invariance to the trajectory of~\eqref{delsys} {\color{black}with initial condition}~$\zeta \mathds{1}$ at time $0$ immediately implies that the vector field at time $0$ evaluated at $\zeta$ must satisfy
\begin{equation}
f(\zeta) + g^0(\zeta \mathds{1}) \le 0. \label{vecfieldineq}
\end{equation}
To make use of~\eqref{vecfieldineq}, we define the {\color{black}functional \mbox{$\bar{g}: \mathcal{C}([-r,0],\mathds{R}^n) \to \mathds{R}^n$},
\begin{equation}
\bar{g}(\phi) := g^0\big(\sup_{\theta \in [-r,0]} \phi(\theta) \mathds{1}\big), \label{eqn:barg}
\end{equation}
{\color{black}where the supremum is taken componentwise, i.e. $\sup_{\theta \in [-r,0]} \phi(\theta)$ is a vector {\color{black}in $\mathds{R}^n$} with components $\sup_{\theta_i \in [-r,0]} \phi_i(\theta_i)$}}. \eqref{eqn:barg} is time-invariant, locally Lipschitz, and {\color{black} satisfies the order-preserving property}~\eqref{gmon}. Thus we can form {\color{black}the}
time-independent monotone~system
\begin{equation}
\frac{dx}{dt}(t) = f(x(t)) + \bar{g}(x_t), \label{boundsys}
\end{equation}
which has an equilibrium at $0$ and is hence positive. Moreover,~\eqref{vecfieldineq} shows that $f(\zeta) + \bar{g}(\zeta \mathds{1}) \le 0$, whence an application of Lemma~\ref{sc522} guarantees that the solution $z(t)$ of~\eqref{boundsys} with initial condition $z_{t_0} = \zeta \mathds{1}$ at any initial time $t_0$ must satisfy\footnote{To be precise, since positivity means that the trajectory is known to be bounded below by $0$ in all components, Lemma~\ref{sc522} gives $\lim_{t\to\infty} z(t) = x^*$ for some $x^*$ satisfying $0 \le x^* \le z(t_0) = \zeta \ll y(0)$. The continuity of the maps $f$ and $\bar{g}$ then means that any such point must satisfy $f(x^*) + \bar{g}(x^* \mathds{1}) = 0$ and so $f(x^*) + g^t(x^* \mathds{1}) = 0$, meaning that it is an equilibrium of the undelayed system~\eqref{undelsys}. But then, because $y(0) \gg 0$ and $\lim_{t \to \infty} y(t) = 0$, the monotonicity of the system~\eqref{undelsys} prohibits this equilibrium from lying in $0 < x^* \le y(0)$, whence we conclude that $x^* = 0$ is the only possibility.} $\lim_{t\to\infty} z(t) = 0$.

We then observe that, by~\eqref{gmon} and the time-invariance of $g^t$ for constant arguments, it must always hold that $g^t(\phi) \le g^t(\sup_{\theta \in [-r,0]} \phi(\theta) \mathds{1}) = g^0(\sup_{\theta \in [-r,0]} \phi(\theta) \mathds{1}) = \bar{g}(\phi)$. Thus, the vector fields of the monotone systems~\eqref{delsys} and~\eqref{boundsys} are related by
\begin{equation} \label{vecineq}
f(\phi(0)) + g^t(\phi) \le f(\phi(0)) + \bar{g}(\phi)
\end{equation}
for all $t \ge 0$ and all $\phi \in \mathcal{C}([-r,0],\mathds{R}^n)$. Inequality~\eqref{vecineq} means that we may invoke Lemma~\ref{st511} to deduce that if $x(t)$ denotes the solution of~\eqref{delsys} through any $x_{t_0} \in \mathcal{C}([-r,0],[0,\zeta])$, then it is guaranteed that $x(t) \le z(t)$ for all $t \ge t_0$. But since $\lim_{t\to\infty}z(t) = 0$ and $x(t)$ must always remain nonnegative, this immediately implies that~$\lim_{t\to\infty}x(t) = 0$. This proves that all trajectories with initial conditions within the set $\mathcal{C}([-r,0],[0,\zeta])$ converge to the origin.
\end{IEEEproof}

We now prove Corollary~\ref{cor1}, which states that when the vector field is continuously differentiable then the region of attraction in Lemma~\ref{theo1} {\color{black}has a nonempty interior and asymptotic stability can be deduced}.

\begin{IEEEproof}[Proof of Corollary~\ref{cor1}]
Since $f$ and $g^t$ are continuously differentiable,~\cite[Remark 3.1.2]{smith} guarantees that $y(t) \gg 0$ for all $t \ge 0$. Let $t^p$ be a time for which $y(t^p) \ll y(0)$. Since $y(t) \gg 0$ for all $t \in [0,t^p]$ and $\lim_{t\to\infty} y(t) = 0$, there must exist $\tau \ge t^p$ such that $T_i(a) = \tau$ for the choice $a = \tfrac{1}{2} \inf_{t\in[0,t^p]} y_i(t) > 0$. This holds for all $i$, so we see that the point~$\zeta$ defined in~\eqref{zdef} satisfies $\zeta \gg 0$.

In addition, stability follows from convergence, using the monotonicity of the system~\eqref{delsys} and the continuous dependence of its solutions on their initial conditions. We thus see that the origin is both stable and attractive on the set $[0,\zeta]$, where $\zeta \gg 0$. We therefore conclude that the origin is an asymptotically stable equilibrium for the delayed system~\eqref{delsys} with arbitrary bounded delays, with a region of attraction containing the nonempty interval~$[0,\zeta]$.
\end{IEEEproof}

{\color{black}Corollary~\ref{cor2} now follows directly from Corollary~\ref{cor1}.}
{\color{black}
\begin{IEEEproof}[Proof of Corollary~\ref{cor2}]
Asymptotically stability of~\eqref{undelsys} means that there exists $\epsilon > 0$ such that for all $y(0) \in (0,\epsilon)$, the solution $y(t)$ of~\eqref{undelsys} satisfies $\lim_{t\to\infty} y(t) = 0$. Corollary~\ref{cor1} thus applies.
\end{IEEEproof}}

We finally show how the method employed in proving Lemma~\ref{theo1} enables us to prove the global stability result stated in Theorem~\ref{theo2}.
Note that continuous differentiability is not required here.

\begin{IEEEproof}[Proof of Theorem~\ref{theo2}]
Consider the modified definitions
\begin{align*}
&T_i(x_i) := \sup \{\tau \colon x_i \le y_i(s) \; \forall s \in (-\infty, \tau]\},\\
&V_i(x_i) := e^{-T_i(x)},\\
&T(x) := \min_i T_i(x_i), \text{ and } V(x) = \max_i V_i(x_i),
\end{align*}
in terms of the solution $y(t)$ of~\eqref{undelsys}. Then $V(x) = e^{-T(x)}$ as before and $V$ is now a well-defined, positive definite, left-continuous function on the whole of $\mathds{R}_+^n$, according~to~the~given~asymptotic properties of $y(t)$ and the continuity of the system~\eqref{undelsys}~\cite{rantzer}. The structure of the proof is then analogous to that of Lemma~\ref{theo1}:

\textbf{Step 1:} As a consequence of the fact that $V$ is now defined on the whole positively invariant orthant $\mathds{R}_+^n$, analogous arguments to those in Step 1 of the proof of Lemma~\ref{theo1} show that all sublevel sets $S_c := \{x \in \mathds{R}_+^n \colon V(x) \le c\}$ are now positively invariant with respect to~\eqref{delsys}.

\textbf{Step 2:} Consider an arbitrary element $x \in S_c$. Following the proof of Lemma~\ref{theo1}, we introduce the mapping $h_i(u) := \inf \{ \tau \ge u \colon \exists \, a \in \mathds{R}_+ \text{ such that } T_i(a) = \tau \}$. Then, for any $x \in S_c$ we have $V(x) \le c$, implying that $T_i(x_i) \ge -\log c$ for all $i$. Thus, it must hold that $T_i(x_i) \ge h_i(-\log c)$, and so $x_i \le y_i(h_i(-\log c))$. Consequently, if we define $\zeta^c \in \mathds{R}_+^n$ by
\begin{equation}
\zeta^c_i := y_i(h_i(-\log c)) \label{zdefthm}
\end{equation}
for all $i$, then $S_c \subseteq [0,\zeta^c]$.

In particular, the sublevel sets are thus bounded, meaning that we can again use the left-continuity and nonincreasing nature of $T_i$ to deduce the existence of $a \in \mathds{R}_+$ such that $h_i(-\log c) = T_i(a)$. Therefore, $a = y_i(h_i(-\log c))$, whence we have all $T_i(\zeta^c_i) = h_i(-\log c) \ge -\log c$ and so $V(\zeta^c) \le c$. This implies that $\zeta^c \in S_c$.

\textbf{Step 3:} As in the proof of Lemma~\ref{theo1}, we therefore have an invariant set $S_c$ containing $\zeta^c$ and contained within $[0,\zeta^c]$, whence we deduce
\begin{equation}
f(\zeta^c) + g^0(\zeta^c \mathds{1}) \le 0. \label{vecfieldineqglobal}
\end{equation}
Furthermore, we observe that since $\lim_{t\to-\infty} y_i(t)~=~\infty$ for all~$i$, the function~$V$ is radially unbounded. Therefore, the sublevel sets satisfy $S_c \to \mathds{R}_+^n$ as $c \to \infty$. In particular, the points $\zeta^c \to \infty$ as $c \to \infty$, meaning that there exist arbitrarily large $\zeta^c$ for which~\eqref{vecfieldineqglobal} holds. Thus, given any $x_{t_0} \in \mathcal{C}([-r,0], \mathds{R}_+^n)$, there exists $\zeta^c$ such that $x_{t_0} \le \zeta^c \mathds{1}$ and~\eqref{vecfieldineqglobal} holds. Then repeating the construction~of~the bounding system~\eqref{boundsys} and applying Lemmas~\ref{st511} and~\ref{sc522}, as in Step 3 of the proof of Lemma~\ref{theo1}, guarantees that the solution $x(t)$ of~\eqref{delsys} with initial condition $x_{t_0}$ must satisfy $\lim_{t\to\infty}x(t) = 0$. Stability now follows by the {\color{black}argument used} in the proof of Corollary~\ref{cor1}. Therefore, the origin is globally asymptotically stable for the delayed system~\eqref{delsys}.\end{IEEEproof}

{\color{black}
\begin{rem}
It should be noted that even though the condition $\lim_{t \to -\infty} y_i(t) = \infty$ for all components $i$ immediately implies the existence of arbitrarily large points at which each component of the vector field is individually nonpositive, it is nontrivial to show that this happens simultaneously for all components. The latter is needed in order to deduce delay-independent stability {\color{black}using Lemma~\ref{sc522} and {\color{black} the monotonicity arguments in Lemma~\ref{st511}}}. It was shown above that the {\color{black}quantity} $\zeta^c$ constructed within the proof of Theorem~\ref{theo2} explicitly defines such points, thus allowing global asymptotic stability to be deduced.
\end{rem}
}

{\color{black}
\section{Discussion} \label{sec:examples}

We have seen already in Theorem~\ref{theo2} how the regions of attraction constructed in Lemma~\ref{theo1} can lead to global convergence results. {\color{black}Here we discuss how our approaches connect with other results in the literature and provide also additional examples where delay-independent stability can be deduced.}

\subsection{Connections with Lemmas~\ref{st511} and~\ref{sc522}} \label{sec:examples2}

Lemmas~\ref{st511} and~\ref{sc522} are important classical results that have been used in a number of interesting studies to prove delay-independent global asymptotic stability for linear, homogeneous, and subhomogeneous monotone systems. These results are based on the fact that in such systems global asymptotic stability of the undelayed system implies the existence of arbitrarily large points where the vector field is negative. Lemma~\ref{theo1} shows that such points can be determined under weaker system assumptions. Even though these points (denoted as~$\zeta$) can in general be difficult to characterize algorithmically, their significance lies in the fact that they allow to prove delay-independent global asymptotic stability for broader classes of systems.

The following corollary and its proof make the connection with Lemmas~\ref{st511} and~\ref{sc522} more explicit. The corollary shows that when a point is known where the vector field is negative, then the delay-independent region of attraction obtained from Lemma~\ref{theo1} has a simple form, recovering arbitrarily closely that obtained from Lemmas~\ref{st511} and~\ref{sc522}. Note also that Lemma~\ref{theo1} additionally allows the delays to be time-varying.

\begin{cor} \label{lemcomprem}
If $v$ is any point in $\mathds{R}_+^n$ such that $f(v) + g^0(v \mathds{1}) \ll 0$, then any solution of the delayed system~\eqref{delsys} with initial condition $x_{t_0} \in \mathcal{C}([-r,0],[0,v))$ satisfies $\lim_{t \to \infty} x(t) = 0$.
\end{cor}

\begin{IEEEproof}
{\color{black}Define} $y(t)$ to be the trajectory of the undelayed system~\eqref{undelsys} through $y(0) = v$. Then the given vector field property implies that, for all sufficiently small $t^p > 0$, $y(t^p) \ll y(0)$ holds and the function defined in~\eqref{hdef} satisfies $h_i(t^p) = t^p$ for all~$i$. Therefore, it follows from~\eqref{zdef} that the point $\zeta$ can be brought arbitrarily close to $v$, and so $[0,v)$ is a region of attraction obtained from Lemma~\ref{theo1}.\end{IEEEproof}

\begin{rem}
{\color{black}Corollary~\ref{lemcomprem}} is also true if the inequality for the vector field is not strict, i.e.  $f(v) + g^0(v \mathds{1}) \leq 0$, and the point $v$ can also be included in the estimate of the region of attraction. This follows directly from Step 3 of the proof of Lemma~\ref{theo1} by replacing $\zeta$ with~$v$.
\end{rem}

\subsection{Subhomogeneous systems} \label{sec:examples3}

We now indicate how the {\color{black}regions} of attraction constructed in Lemma~\ref{theo1} can also recover known delay-independent stability results for monotone systems that are also subhomogeneous. To see this, let us first recall the definition of subhomogeneity, as used in~\cite{bokharaie2}.

\begin{defn} \normalfont
The functions $f$ and $g^t$ are said to be subhomogeneous of degree $\alpha > 0$ if they satisfy $f(\lambda x) \le \lambda^\alpha f(x)$ and $g^t(\lambda \phi) \le \lambda^\alpha g^t(\phi)$ for all $\lambda \ge 1$, all $t \ge 0$, all $x \in \mathds{R}_+^n$, and all $\phi \in \mathcal{C}([-r,0],\mathds{R}_+^n)$.
\end{defn}

\begin{cor} \label{homogcompcor}
Suppose that $f$ and $g^t$ are both continuously differentiable for all fixed $t$ and subhomogeneous of any degree $\alpha > 0$ and that the origin is known to be globally asymptotically stable for the undelayed system~\eqref{undelsys}. Then the origin is globally asymptotically stable for the delayed system~\eqref{delsys}.
\end{cor}
\begin{IEEEproof}
Consider the time-invariant bounding system~\eqref{boundsys}. This is then positive, monotone, and subhomogeneous, and the origin of the corresponding undelayed system is known to be globally asymptotically stable, whence we recall from the proof of Theorem~3.2 in~\cite{bokharaie2} that given any $x_0 \ge 0$, there exists $v \gg x_0$ such that $f(v) + \bar{g}(v \mathds{1}) \ll 0$. Therefore, it follows by definition {\color{black}of $\bar g$ given in \eqref{eqn:barg}} that $f(v) + g^0(v \mathds{1}) \ll 0$, and so Corollary~\ref{lemcomprem} guarantees that the region of attraction of the origin of the delayed system~\eqref{delsys} contains the region~$[0,v)$. Allowing $x_0 \to \infty$, $v$ can then be made arbitrarily large, meaning that the convergence region obtained from Lemma~\ref{theo1} is the entire positive orthant $\mathds{R}_+^n$.
\end{IEEEproof}

{\color{black}{\color{black}The proofs of Corollaries~\ref{lemcomprem} and \ref{homogcompcor} show that} subhomogeneous monotone systems admit undelayed trajectories through which Lemma \ref{theo1}} can be used to deduce delay-independent global asymptotically stability (with time-varying delays), recovering Corollary~3.3 from~\cite{bokharaie2}. It should be noted, though, as was shown in~\cite{bokharaie2}, that subhomogeneity is a sufficiently strong property that leads to delay-independent stability also for classes of systems that are not~monotone.

\subsection{Global asymptotic stability {\color{black}via Gronwall's Lemma}} \label{sec:examplesgas}

{\color{black}We now give an example where the regions of attraction quantified in Lemma~\ref{theo1} can be used to deduce delay-independent global asymptotic stability {\color{black}in} a system where Theorem~\ref{theo2} does not apply\footnote{\color{black}Note that this example also does not satisfy~\cite[Assumption A]{bokharaiecorrec}, which was used in~\cite{bokharaiecorrec} to deduce global asymptotic stability for planar systems. Additionally, because both components in~\eqref{ex3system} cannot simultaneously tend to~$\infty$ in negative time, Theorem~\ref{theo2} also cannot be applied for this system.} and {\color{black}which} is also not subhomogeneous.} In particular we consider a delayed form of the second example in~\cite{rantzer}, with delays $\tau_1(t),\tau_2(t) \in [0,r]$,
\begin{equation} \label{ex3system}
\begin{aligned}
&\frac{d x_1}{dt}(t) = -x_1(t) + \frac{x_2(t-\tau_1(t))^2}{x_2(t-\tau_1(t))^2+1} \\
&\frac{d x_2}{dt}(t) = x_1(t-\tau_2(t)) - \frac{2x_2(t)^2}{x_2(t)^2+1}.
\end{aligned}
\end{equation}
{\color{black}The equations \eqref{ex3system} define a {\color{black}monotone system for initial conditions taking values in $\mathds{R}_+^n$}, and} the origin of the corresponding undelayed system is globally asymptotically stable\footnote{Recall that the stability definitions considered are restricted to the nonnegative orthant $\mathds{R}_+^n$.} as discussed in~\cite{rantzer}. {\color{black}For any {\color{black}$\kappa \ge 1$}, consider the initial~condition $y(0) = (\tfrac{3\kappa^2-1}{k^2+1},\kappa)^T$ in Lemma~\ref{theo1}. It can be shown\footnote{\color{black}Consider the undelayed system corresponding to \eqref{ex3system}. Whenever $x_1 \ge 1$, we have $\tfrac{dx_1}{dt} < 0$. Also, $\tfrac{dx_2}{dt} > 0$ when $x_1 \in \big(\tfrac{2x_2^2}{x_2^2+1}, \tfrac{3x_2^2-1}{x_2^2+1}\big)$ and $\tfrac{dx_2}{dt} < 0$ when $x_1 \in \big(1, \tfrac{2x_2^2}{x_2^2+1}\big)$. Furthermore, for any $\eta \in \big(0,\tfrac{2\xi^2}{\xi+1} - 1\big)$ and $\tilde{\xi} \ge \xi \ge 1$, it holds that ${\left|\tfrac{dx_2}{dx_1}\right|_{ x_1 = \tfrac{2\xi^2}{\xi^2+1} + \eta \text{, } x_2 = \xi} < \left|\tfrac{dx_2}{dx_1}\right|_{x_1 = \tfrac{2\xi^2}{\xi^2+1} - \eta \text{, } x_2 = \tilde{\xi}}}$. The above imply, for the trajectory $y(t)$ through the given initial conditions, that {\color{black}the initial increase in $y_2$ until $y_1=2y_2^2/(y_2^2+1)$ is smaller than its subsequent decrease until $y_1=1$. Therefore,} there exists $t^p \ge 0$ for which $(1,\kappa-1)^T \ll y(t^p) \ll y(0)$, whence the point $\zeta$ defined in \eqref{zdef} can be chosen such \color{black} that $\zeta \gg (1,\kappa-1)^T$.} } that the points $\zeta$ can be chosen to satisfy $\zeta \gg (1,\kappa-1)^T$. Therefore, by {\color{black}choosing $\kappa$ arbitrarily large}, Lemma~\ref{theo1} gives delay-independent convergence on a region of attraction containing the set $\{x \colon x_1 \le 1\}$. {\color{black}{\color{black}Additionally}, it can be deduced from an application of Gronwall's Lemma that all trajectories of the delayed system must reach this set in finite {\color{black}time}\footnote{To see this, observe that whenever $x_1(t) \ge \frac{3}{2}$, we have $\frac{dx_1}{dt}(t) < -\frac{1}{2}$, implying that it takes only finite time to reach the set $\{x \colon x_1 \le \frac{3}{2}\}$. For any given trajectory of the delayed system, call this time~$\tau$. Thereafter, the trajectory always remains within this set and hence $\frac{dx_2}{dt}(t) \le \frac{3}{2} - \frac{2x_2(t)^2}{x_2(t)^2+1}$, which implies that $x_2(t) < b$ for all $t \ge \tau$, for some bounding constant~$b$. Thus, we have $x_2(t) \le B := \max \{ \sup_{t \in [0,\tau]} x_2(t), b \} < \infty$. Consequently, $\frac{dx_1}{dt}(t) \le -x_1(t) + \frac{B^2}{B^2+1}$ for all $t \ge 0$ and therefore, since $\frac{B^2}{B^2+1} < 1$, Gronwall's Lemma establishes the desired property.}.} Consequently, all solutions eventually enter a region on which Lemma~\ref{theo1} guarantees convergence. {\color{black}Therefore, the origin in the system~\eqref{ex3system} is globally asymptotically stable for arbitrary bounded time-varying delays.}
}

\section{Conclusions} \label{sec:conclusions}

Within this paper we have considered the question of when it is possible to guarantee delay-independent stability for nonlinear monotone systems of delay differential equations. We showed that, when the undelayed system admits a convergent solution that is unbounded in all components in negative time, then this is sufficient to deduce that the corresponding delayed system will be globally asymptotically stable under arbitrary bounded delays, which can be heterogeneous and time-varying. This result followed from a more general result in which we showed {\color{black}that positive convergent trajectories of the undelayed system explicitly determine delay-independent regions} on which all solutions of the corresponding delayed system must be convergent. These results demonstrate that it is possible to use the quasimonotonicity property alone to infer information about the behavior of a system under arbitrary bounded delays directly from knowledge of the {\color{black}asymptotic behavior of} the corresponding undelayed system.

\section*{\color{black}Appendix A\\Quasimonotonicity and comparison properties}

The monograph~\cite{smith} provides an excellent reference on the general theory of monotone systems. Chapter $5$ focuses on monotone systems of delay differential equations, analyzing their properties in detail and eventually proving a powerful result of generic convergence for autonomous systems under assumptions of boundedness, continuous differentiability, and irreducibility. Several of the preliminary results, however, are stated for general monotone delayed systems, making them useful here. We will now review two useful comparison results that are invoked within the proofs of our main results in Section~\ref{sec:proofs}.

We shall consider the general delay differential equation
\begin{equation}
\frac{dx}{dt}(t) = f(t,x_t), \label{appsys}
\end{equation}
with $f: \Lambda \times \mathcal{C}([-r,0],\Omega) \to \mathds{R}^n$ continuous on open subsets $\Lambda \subseteq \mathds{R}$ and $\Omega \subseteq \mathds{R}^n$. We say that $f$ satisfies the quasimonotonicity property~if
\begin{equation}
\phi \le \psi \text{ and } \phi_i(0) = \psi_i(0) \Rightarrow f_i(t,\phi) \le f_i(t,\psi) \label{appmon}
\end{equation}
for all $t \in \Lambda$ and all $\phi, \psi \in \mathcal{C}([-r,0],\Omega)$. Given $t_0 \in \Lambda$ and $\phi \in \mathcal{C}([-r,0],\Omega)$, we let $x(t,t_0,\phi,f)$ denote the maximally defined solution that satisfies~\eqref{appsys} for all $t \ge t_0$ and passes through $x_{t_0} = \phi$.

{\color{black}It is well-known that systems that satisfy the property of quasimonotonicity are monotone, i.e. if two initial conditions $\phi, \psi \in \mathcal{C}([-r,0],\Omega)$ satisfy $\phi \le \psi$ then $x(t,t_0,\phi,f) \le x(t,t_0,\psi,f)$ for $t \ge t_0$. This is stated in the following lemma (Theorem 5.1.1 in~\cite{smith}) which is the first result we recall. Note that this lemma also allows trajectories of systems with different vector fields to be compared.}

\begin{lem} \label{st511}
Let $f,g : \Lambda \times \mathcal{C}([-r,0],\Omega) \to \mathds{R}^n$ be continuous in $t$ and locally Lipschitz in $\phi$, and assume that either $f$ or $g$ satisfies~\eqref{appmon}. Assume also that $f(t,\phi) \le g(t,\phi)$ for all $t \in \Lambda$ and all $\phi \in \mathcal{C}([-r,0],\Omega)$. If $t_0 \in \Lambda$ and $\phi, \psi \in \mathcal{C}([-r,0],\Omega)$ satisfy $\phi \le \psi$, then $x(t,t_0,\phi,f) \le x(t,t_0,\psi,g)$ holds for all $t \ge t_0$ for which both are defined.
\end{lem}

The second result guarantees monotonic convergence of any bounded trajectory of an autonomous system at whose initial condition the vector field is either nonnegative or nonpositive in all components. This is Corollary 5.2.2 in~\cite{smith}.

\begin{lem} \label{sc522}
Let $f : \mathcal{C}([-r,0],\Omega) \to \mathds{R}^n$ be time-invariant, {\color{black}locally} Lipschitz, and satisfy~\eqref{appmon}. If $v \in \Omega$ is such that $f(v \mathds{1}) \ge 0$ ($f(v \mathds{1}) \le 0$), then $x(t,t_0,v \mathds{1},f)$ is nondecreasing (nonincreasing) in $t \ge t_0$. If the positive orbit of $h \mathds{1}$ has compact closure in $\Omega$, then there exists $k \ge v$ ($k \le v$) such that $x(t,t_0,v\mathds{1},f) \to k$ as $t \to \infty$.
\end{lem}

{\color{black}\section*{Appendix B\\Nonpositive systems}}

The approach within the main paper focused, for simplicity of presentation, on monotone systems with an equilibrium at the origin. We now demonstrate that the analysis in Section \ref{sec:convres} can also be applied to yield delay-independent convergence guarantees in general monotone systems of the form of~\eqref{undelsys} and~\eqref{delsys} admitting an arbitrary equilibrium point $x^* \in \mathbb{R}^n$. We thus use within this section the standard definitions for stability as stated, for example, in Definition 5.1.1 in \cite{hale}.

The monotonicity properties~\eqref{fmon} and~\eqref{gmon} imply that the orthants $x^* + \mathbb{R}_+^n$ and $x^* - \mathbb{R}_+^n$ must be positively invariant for both systems~\eqref{undelsys} and~\eqref{delsys}. It is therefore natural that coordinate transformations can enable Theorem~\ref{theo2} to be extended to the following more general result.

\begin{cor} \label{nonposcor1}
Suppose that the undelayed system~\eqref{undelsys} admits solutions $\overline{y}(t)$ and $\underline{y}(t)$ satisfying $\lim_{t \to -\infty} \overline{y}_i(t)$ $= \infty$, $\lim_{t \to \infty} \overline{y}_i(t) = x^*$, $\lim_{t \to -\infty} \underline{y}_i(t) = -\infty$, and $\lim_{t \to \infty} \underline{y}_i(t) = x^*$ for all components $i$. Then $x^*$ is globally asymptotically stable for the delayed system~\eqref{delsys}.
\end{cor}

Analogously to Section \ref{sec:convres}, this represents a limiting case of a more general technical result which extends Lemma \ref{theo1} to quantify delay-independent regions of attraction. Note that the points $\overline{\zeta},\underline{\zeta}$ mentioned in the corollary are defined within the proof of the corollary in a way analogous to the point $\zeta$ in  Lemma \ref{theo1}.

\begin{cor} \label{nonposcor2}
Suppose that the undelayed system~\eqref{undelsys} admits solutions $\overline{y}(t)$ and $\underline{y}(t)$ satisfying $\overline{y}_i(0) > x^*$, $\lim_{t \to \infty} \overline{y}_i(t) = x^*$, $\underline{y}_i(0) < x^*$, and $\lim_{t \to \infty} \underline{y}_i(t) = x^*$ for all components $i$. Then there exist points $\overline{\zeta} \ge x^*$ and $\underline{\zeta} \le x^*$, defined in terms of the trajectories $\overline{y}(t)$ and $\underline{y}(t)$ respectively, such that any solution of the delayed system~\eqref{delsys} with initial condition $x_{t_0} \in \mathcal{C}([-r,0],[\underline{\zeta},\overline{\zeta}])$ satisfies $\lim_{t \to \infty} x(t)~=~x^*$.
\end{cor}

If the vector field is additionally continuously differentiable, Corollary \ref{nonposcor2} again extends as follows to enable stronger conclusions to be drawn.

\begin{cor} \label{poscor1}
Suppose that $f$ and $g^t$ are both continuously differentiable for all fixed $t$ {\color{black}and that the undelayed system~\eqref{undelsys} admits solutions $\overline{y}(t)$ and $\underline{y}(t)$ satisfying $\overline{y}_i(0) > x^*$, $\lim_{t \to \infty} \overline{y}_i(t) = x^*$, $\underline{y}_i(0) < x^*$, and $\lim_{t \to \infty} \underline{y}_i(t) = x^*$ for all components $i$}. Then the points  $\overline{\zeta},\underline{\zeta}$ in Corollary~\ref{nonposcor2} can be chosen such that $\overline{\zeta} \gg x^*$ and $\underline{\zeta} \ll x^*$, and the point $ x^*$ is asymptotically stable for the delayed system~\eqref{delsys}.
\end{cor}

Corollary~\ref{poscor1} then also immediately implies the following simple corollary.

\begin{cor} \label{poscor2}
Suppose that $f$ and $g^t$ are both continuously differentiable for all fixed $t$ and that the point $ x^*$ is asymptotically stable for the undelayed system~\eqref{undelsys}. Then $ x^*$ is also asymptotically stable for the delayed system~\eqref{delsys}.
\end{cor}

\balance

{\color{black}We now show how the results stated above can be derived from the theory established in Sections \ref{sec:convres} and \ref{sec:proofs}.}
\begin{IEEEproof}[Proof of Corollary \ref{nonposcor2}.]
Make the coordinate change $u = x-x^*$ in~\eqref{undelsys} and~\eqref{delsys}. This gives
\begin{align}
\frac{du}{dt}(t) &= f(x^*+u(t)) + g^t(x^* \mathds{1}+u(t) \mathds{1}) \nonumber \\
&=\tilde{f}(u(t)) + \tilde{g}^t(u(t)\mathds{1}) \label{undelsysmod}
\end{align}
and
\begin{align}
\frac{du}{dt}(t) &= f(x^*+u(t)) + g^t(x^* \mathds{1}+u_t) \nonumber \\
&= \tilde{f}(u(t)) + \tilde{g}^t(u_t). \label{delsysmod}
\end{align}
The systems~\eqref{undelsysmod} and~\eqref{delsysmod} then have an equilibrium at the origin and, moreover, their nonlinearities $\tilde{f}$ and $\tilde{g}^t$ clearly satisfy the monotonicity properties~\eqref{fmon} and~\eqref{gmon}. Finally, using the given trajectory $\overline{y}(t)$ of~\eqref{undelsys}, we see that~\eqref{undelsysmod} admits a solution $\overline{v}(t)$ defined by $\overline{v}(t) = \overline{y}(t) - x^*$ for all $t \ge 0$ satisfying $\overline{v}_i(0) > 0$ and $\lim_{t \to \infty} \overline{v}_i(t) = 0$ for all components $i$. The analysis of Lemma~\ref{theo1} can thus be applied, yielding a point $\overline{\gamma} \ge 0$ defined as in \eqref{zdef} such that the origin of~\eqref{delsysmod} is asymptotically stable with region of attraction containing~$[0,\overline{\gamma}]$.

Analogous analysis with the positive monotone systems
\begin{align}
\frac{du}{dt}(t) &= -f(x^*-u(t)) - g^t(x^* \mathds{1}-u(t) \mathds{1}) \nonumber \\
&=\tilde{f}(u(t)) + \tilde{g}^t(u(t)\mathds{1}) \label{undelsysmod2}
\end{align}
and
\begin{align}
\frac{du}{dt}(t) &= -f(x^*-u(t)) - g^t(x^* \mathds{1}-u_t) \nonumber \\
&= \tilde{f}(u(t)) + \tilde{g}^t(u_t), \label{delsysmod2}
\end{align}
obtained under the coordinate change $u = x^*-x$, gives a point $\underline{\gamma} \ge 0$ such that~\eqref{delsysmod2} has asymptotically stable origin on the region of attraction $[0,\underline{\gamma}]$.

Combining these analyses and denoting $\overline{\zeta} =  x^*+\overline{\gamma}$ and $\underline{\zeta} =  x^*-\underline{\gamma}$ shows that the region of attraction of the equilibrium point $x^*$ of the original delayed system~\eqref{delsys} contains the union $[\underline{\zeta},x^*] \cup [x^*,\overline{\zeta}]$. Finally, let $\overline{z}(t)$ and $\underline{z}(t)$ respectively denote the trajectories of \eqref{delsys} with initial conditions $\overline{\zeta} \mathds{1}$ and $\underline{\zeta} \mathds{1}$ and let also $\underline{\zeta} \le x_{t_0} \le \overline{\zeta}$. Then the fact that $\overline{z}(t)$ and $\underline{z}(t)$ converge to $x^\ast$ in conjunction with Lemma~\ref{st511} guarantees that $\lim_{t \to \infty} x(t) = x^*$ also. Therefore, the region of attraction must in fact contain~$[\underline{\zeta},\overline{\zeta}]$. 
\end{IEEEproof}

\begin{IEEEproof}[Proof of Corollary \ref{poscor1}]
Analogously to the proof of Corollary~\ref{cor1}, we invoke~\cite[Remark 3.1.2]{smith} to deduce that the points in Corollary \ref{nonposcor2} satisfy $\overline{\zeta} \gg x^* \gg \underline{\zeta}$. 
\end{IEEEproof}

{\color{black}\begin{IEEEproof}[Proof of Corollary \ref{poscor2}]
As in the proof of Corollary~\ref{cor2}, the asymptotic stability of \eqref{undelsys} implies the existence of $\overline{\epsilon} > 0$ such that $\lim_{t\to\infty} y(t) = 0$ for any solution $y(t)$ of \eqref{undelsys} with initial condition $y(0) \in (x^*,x^*+\overline{\epsilon})$. Analogously, there exists $\underline{\epsilon} > 0$ such that $\lim_{t\to\infty} y(t) = 0$ whenever $y(0) \in (x^*-\overline{\epsilon},x^*)$. Corollary \ref{poscor1} can then be applied to deduce the result.
\end{IEEEproof}}

\begin{IEEEproof}[Proof of Corollary \ref{nonposcor1}.]
Analogously to the proof of Corollary \ref{nonposcor2}, the trajectories $\overline{v}(t) = \overline{y}(t) - x^*$ and $\underline{v}(t) = x^* - \underline{y}(t)$ satisfy the asymptotic properties required in Theorem \ref{theo2} for the positive monotone undelayed systems \eqref{undelsysmod} and \eqref{undelsysmod2}, respectively. Therefore, applying Theorem~\ref{theo2} and invoking Lemma \ref{st511} as before immediately implies the desired result. 
\end{IEEEproof}

\end{document}